\documentclass[12pt,a4paper,reqno]{amsart}
\usepackage{amsmath}
\usepackage{geometry}
\usepackage{amssymb, latexsym}
\usepackage{verbatim}
\usepackage{enumerate}
\usepackage{mathrsfs}  
\usepackage{comment}
\usepackage{txfonts}

\usepackage{mathtools}
\usepackage[tableposition=top]{caption}
\usepackage{booktabs,dcolumn}

\theoremstyle{plain}

\newtheorem{theorem}{Theorem}[section]

\theoremstyle{definition}

\newtheorem{definition}[theorem]{Definition}
\newtheorem{remark}[theorem]{Remark}

\newtheorem{example}[theorem]{Example}

\begin{document}

\title[]{Sheaves and conditional sets}

\author{Asgar Jamneshan}
\address{Departement of Mathematics, University of California, Los Angeles.}
\email{jasgar@math.ucla.edu}

\begin{abstract} 
The purpose of this note is to record a connection between sheaves on complete Boolean algebras and conditional sets.  
This connection yields a transfer principle for conditional set theory. 
On the other hand we use conditional set theory to interpret the internal objects of a respective sheaf topos. 
\end{abstract}

\maketitle

\setcounter{tocdepth}{1}



%


\subsection{Sheaves on Boolean algebras}\label{s11}
We collect basic facts about sheaves on Boolean algebras; relying on  \cite{maclane2012sheaves}.  
Details with an exact referencing to the relevant parts of \cite{maclane2012sheaves} are worked out in the thesis of the author in \cite[Chapter 6.1]{jamneshan2014theory}.  
Fix a complete Boolean algebra\footnote{For any unexplained definition or result in the theory of Boolean algebras and measure algebras, the interested reader is referred to  \cite{monk1989handbook}.}   $\mathcal{A}=(\mathcal{A},\wedge,\vee,{}^c,0,1)$. 
To define a sheaf on $\mathcal{A}$ by following the abstract constructions in \cite{maclane2012sheaves}, it is necessary to view $\mathcal{A}$ as a category. 
To this end, consider as objects the elements of $\mathcal{A}$ and an arrow exists between two objects $a$ and $b$ whenever $a\leq b$ where this partial order is defined by $a\wedge b=a$; let $\Delta$ denote the set of all arrows of $\mathcal{A}$. 
Then $\mathcal{A}$ is a Cartesian closed category 
with initial and terminal objects $0$ and $1$. 
Moreover, the pullback of $(a,c)$ and $(b,c)$ is $a\wedge b$. 
Therefore, a \emph{sieve} on $a\in\mathcal{A}$ is a family $\{(a_i,a)\colon i\in I\}\subseteq \Delta$ such that $\{a_i\colon i\in I\}$ is downwards closed, that is if $b\in\mathcal{A}$ with $b\leq a_i$ for some $i$, then $b\in \{a_i\colon i\in I\}$.  
The \emph{Grothendieck sup topology} on $\mathcal{A}$ is the function $J$ which assigns to each $a\in\mathcal{A}$ the set 
\[
J(a):=\left\{ \{ (a_i,a)\colon i\in I\}\text{ sieve such that} \vee_{i\in I} a_i=a\right\}. 
\]
Let $P$ assign to each $a\in\mathcal{A}$ the set of partitions of $a$: 
\[
P(a):=\{ \{(a_i,a)\colon i\in I\} \subseteq \Delta \colon a_i\wedge a_j=0 \text{ if } i\neq j \text{ and}  \vee a_i=a\}. 
\]
Note that $0$ is allowed to be an element of a partition.  
Then $P$ is a Grothendieck basis on $\mathcal{A}$ which generates the sup topology $J$. 
Indeed, let $a\in\mathcal{A}$ and $\{(a_i,a):i\in I\}$ be a sieve with $\vee_{i\in I} a_i=a$. 
Assume a well-ordering on the index $I$. 
Then $b_i:=a_i\wedge \left(\vee_{j<i} b_j\right)^c$, $i\in I$, defines a family in $P(a)$ such that $\{(b_i,a)\colon i\in I\}\subseteq \{(a_i,a)\colon i\in I\}$.  
Therefore, any sheaf on the site $(\mathcal{A},J)$ can be associated uniquely to a sheaf on the site $(\mathcal{A},P)$. 
For $\{(a_i,a)\colon i\in I\}\in P(a)$, we write $(a_i)=(a_i)_{i\in I}$ henceforth.  
\begin{definition}\label{def:sheaf}
A pair $\mathcal{X}=((X_a)_{a\in\mathcal{A}}, ( f^b_a)_{(a,b)\in \Delta})$ consisting of a family of sets $X_a$, $a\in\mathcal{A}$, and a family of functions $f^b_a\colon X_b\to X_a$, $(a,b)\in\Delta$, is called a \emph{sheaf} on $(\mathcal{A},P)$ if 
\begin{itemize}
\item[(i)] $ f^a_a$ is the identity function for every $a\in\mathcal{A}$;
\item[(ii)] $ f^b_a \circ f^c_b = f^c_a$ for all $a,b,c\in\mathcal{A}$ with $a\leq b\leq c$; 
\item[(iii)] for each $a\in\mathcal{A}$, every $(a_i)_{i\in I}\in P(a)$, and all families $(x_i)_{i\in I}\in \prod_{i\in I} X_{a_i}$ there exists a unique element $x\in X_a$ such that $ f^a_{a_i}(x)=x_i$ for all $i\in I$. This unique element $x$ is called  the \emph{amalgamation} of the family $(x_i)_{i\in I}$, and the notation $x=\sum_{i\in I} x_i|a_i=\sum x_i|a_i$ will be used. 
\end{itemize}
Let $\textbf{Sh}(\mathcal{A},P)$ denote the Grothendieck topos of all sheaves on the site $(\mathcal{A},P)$. 
\end{definition}
The following simple properties of a sheaf will be useful later. 
\begin{itemize}
\item[(a)] $X_0$ is a singleton. This follows from the conventions $\vee_\emptyset=0$ and $\prod_\emptyset=\{\ast\}$ and the uniqueness of amalgamations. 
\item[(b)] There exists a largest $a_\ast\in\mathcal{A}$ such that $X_{a_\ast}\neq \emptyset$.  
Indeed, let $M=\{a\in\mathcal{A}\colon X_a\neq \emptyset\}$.  
If $X_b\neq \emptyset$ and $a\leq b$, then $X_a\neq \emptyset$. 
Thus, one finds $(a_i)\in P(\vee M)$ such that $a_i\leq a$ for some $a\in M$ for all $i$, and the claim follows from existence of amalgamations. 
Denote by $a_\mathcal{X}:=a_\ast$. 
\item[(c)] If $(a,b)\in \Delta$ and $X_b \neq\emptyset$, then $ f^b_a$ is surjective. We may assume $a< b$. Pick $y\in X_a$ and any $z\in X_{a^c\wedge b}$. Then there exists a unique $x\in X_b$ such that $ f^b_a(x)=y$ and $f^b_{a^c\wedge b}(x)=z$.   
\end{itemize}
\begin{definition}
Let  $\mathcal{X}=((X_a)_{a\in\mathcal{A}}, (f^b_a)_{(a,b)\in \Delta})$ and $\mathcal{Y}=((Y_a)_{a\in\mathcal{A}}, (g^b_a)_{(a,b)\in \Delta})$ be a sheaves, and  $\mathcal{X}_i=((X_a(i))_{a\in\mathcal{A}}, (f^b_a(i))_{(a,b)\in \Delta})$, $i\in I$, be a family of sheaves. 
\begin{itemize}
\item For $c\in\mathcal{A}$, let $\mathcal{A}_c:=\{a\in\mathcal{A}\colon a\leq c\}$ be a relative algebra of $\mathcal{A}$, and denote by $\Delta_c=\{(a,b)\in \Delta\colon a\leq b\leq c\}$. 
Then $\mathcal{X}|c:=((X_a)_{a\in\mathcal{A}_c}, (f^b_a)_{(a,b)\in \Delta_c})$ is a sheaf on the site $(\mathcal{A}_c,P)$, and is called the \emph{restriction} of $\mathcal{X}$ to $\mathcal{A}_c$. 
\item A sheaf $\mathcal{Y}=((Y_a)_{a\in\mathcal{A}},(g^b_a)_{(a,b)\in\Delta})$ is a \emph{subsheaf} of $\mathcal{X}$ if 
\begin{itemize}
\item[(i)] $Y_a\subseteq X_a$ for all $a\in\mathcal{A}$;
\item[(ii)] $g^b_a$ is the restriction of $ f^b_a$ to $Y_b$ for all $(a,b)\in \Delta$;
\item[(iii)] if for $a\in\mathcal{A}$, $(a_i)_{i\in I}\in P(a)$ and $x\in X_a$ one has $ f^a_{a_i}(x)\in Y_{a_i}$ for all $i\in I$, then $x\in Y_a$.  
\end{itemize}
For $a\in\mathcal{A}$, define 
\[
\mathcal{P}_a:=\{\mathcal{Y}|a\colon \mathcal{Y} \text{ subsheaf of }\mathcal{X}\}. 
\]
For $(a,b)\in \Delta$, define 
\[
\gamma^b_a\colon \mathcal{P}_b\to\mathcal{P}_a, \quad \gamma^b_a(\mathcal{Y}|b):=\mathcal{Y}|a. 
\]
Then $\mathfrak{P}(\mathcal{X})=((\mathcal{P}_a)_{a\in\mathcal{A}},(\gamma^b_a)_{(a,b)\in\Delta})$ is a sheaf on $(\mathcal{A},P)$, and the \emph{power object} of $\mathcal{X}$ (while $\mathcal{A}$ is the \emph{subobject classifier}).  
On $\mathcal{P}(\mathcal{X})$ the \emph{subobject relation} is defined by 
\[
\mathcal{Y}\leq \mathcal{Z} \text{ if and only if } Y_a\subseteq Z_a \text{ for all }a\in\mathcal{A}. 
\]
\item The \emph{Cartesian product} $\prod_{i\in I} \mathcal{X}_i$ consists of the family of sets 
\[
\prod_{i\in I} X_a(i), \quad a\in\mathcal{A}, 
\]
and of the functions 
\[
\alpha\colon \prod_{i\in I} X_b(i) \to \prod_{i\in I} X_a(i), \quad \alpha((x(i))_{i\in I}):= (f^b_a(i)(x_i))_{i\in I}
\]
\item A \emph{sheaf morphism} from $\mathcal{X}$ to $\mathcal{Y}$ is a family of functions $F_a:X_a\to Y_a$, $a\in\mathcal{A}$, such that 
\[
g^b_a\circ F_b= F_a \circ f^b_a\quad \forall (a,b)\in \Delta. 
\]
Two sheaves $\mathcal{X}$ and $\mathcal{Y}$ are said to be \emph{isomorphic} if there exists a morphism $F$ from $\mathcal{X}$ to $\mathcal{Y}$ each component of which is bijective.
There is no sheaf morphism if $a_\mathcal{X}\wedge b_\mathcal{Y}^c>0$ since there does not exist a function from a nonempty set into the empty set.  
Otherwise, if $a_\mathcal{X}\leq b_\mathcal{Y}$, then the component $F_a$ is the empty function for all $a\in\mathcal{A}$ with $a\wedge a_\mathcal{X}^c>0$.  
Suppose $a\leq a_\mathcal{X}\leq b_\mathcal{Y}$.  Then for all $(a_i)\in P(a)$ and $(x_i)_{i\in I}\in \prod_{i\in I} X_{a_i}$ one has 
\[
F_a\left(\sum x_i|a_i\right)=\sum F_a(x_i)|a_i
\]
\item Any classical set $X$ can be \emph{sheafified}. For $a\in\mathcal{A}$, let $X_a$ denote the set of all families $\{(x_i,a_i)\colon i\in I\}$ such that two such families $\{(x_i,a_i)\colon i\in I\}$ and $\{(y_j,b_j)\colon j\in J\}$ are equivalent if $\vee\{a_i\colon x_i=x\}=\vee\{b_j\colon y_j=x\}$  for all $x\in X$. Define $f^b_a\colon X_b\to X_a$ by $f^b_a(\{(x_i,b_i)\colon i\in I\})=\{(x_i,b_i\wedge a)\colon i\in I\}$ for each $(a,b)\in \Delta$.  
For example, the sheafification of the classical natural numbers is the natural numbers object in $\textbf{Sh}(\mathcal{A},P)$.  
\end{itemize}
\end{definition}
As an important consequence of the following theorem we obtain an interpretation of ``bounded Zermelo'' set theory in the \emph{internal logic} of  $\textbf{Sh}(\mathcal{A},P)$ while truth is evaluated in $\mathcal{A}$ (viewed as a sheaf on itself), see e.g. \cite{shulman2019comparing} for a recent discussion.  
\begin{theorem}
The Grothendieck topos $\textbf{Sh}(\mathcal{A},P)$ is Boolean, has a natural numbers object, and satisfies the axiom of choice. 
Moreover,  $\textbf{Sh}(\mathcal{A},P)$ is well-pointed iff $\mathcal{A}=\{0,1\}$ in which case it collapses to the category of sets. 
\end{theorem}

\subsection{Conditional sets}\label{s12}

We give an introduction into conditional sets which is necessary to prove an equivalence to the topos of sheaves $\textbf{Sh}(\mathcal{A},P)$, details and proofs can be found in \cite{drapeau2016algebra} (while keeping in mind Remark \ref{superfluous}). 
 
\begin{definition}
A \emph{conditional set} is a collection $\mathbf{X}$ of objects $x|a$, where $x$ is an element of a  nonempty  set $X$ and $a\in\mathcal{A}$,  with the following properties:   
    \begin{itemize}
        \item[(i)] If $x|b=y|b$ and $a \leq b$ then $x|a=y|a$;
        \item[(ii)] if $(a_i)_{i\in I}\in P(1)$ and $(x_i)_{i\in I}$ is a family in $X$, then there exists exactly one $x\in X$ such that $x|a_i=x_i|a_i$ for all $i$.  
        This unique $x$ is called the concatenation of $(x_i)$ (with respect to $(a_i)$), and is denoted by $x=\sum_{i\in I} x_i|a_i=\sum x_i|a_i$. 
	\end{itemize}
	The symbol $|$ is called a restriction operation. The classical  nonempty  set over which the structure of a conditional set is defined will be called its carrier.  
A  nonempty  subset $Y\subseteq X$ is said to be stable under concatenations, or \emph{stable} for short, if it contains all concatenations of its elements with respect to any partition of unity.  
A stable set $Y$ gives rise to a new conditional set $\{y|a\colon y\in Y,\; a\in\mathcal{A}\}$ where the restriction operation is inherited from $X$.  
For $a\in\mathcal{A}$, let $X_a:=\{x|a\colon x\in X\}$, and define 
\[
\mathbf{X}|a:=\{(x|a)|b\colon x|a\in X_a, \; b\leq a\}
\]
where $(x|a)|b:=x|b$.  By \cite[Proposition 2.4]{drapeau2016algebra}, the collection $\mathbf{X}|a$ is a conditional set on the relative algebra $\mathcal{A}_a=\{b\in\mathcal{A}\colon b\leq a\}$, named the \emph{restriction} of $\mathbf{X}$ to $a$.  
Therefore, for any stable subset $Y\subseteq X$ and $a\in \mathcal{A}$, one can define the conditional set $\mathbf{Y}|a$ on $\mathcal{A}_a$, and any conditional set of this form is called a \emph{conditional subset} of $\mathbf{X}$.  

Let $\mathcal{S}(\mathbf{X})$ denote the collection of stable subsets of $X$. We may endow the classical set $\mathcal{S}(\mathbf{X})$ with the structure of a conditional set via the restriction operation $Y|a:=\{y|a\colon y\in Y\}$. By \cite[Proposition 2.4]{drapeau2016algebra}, the unique concatenation of a family $(Y_i)_{i\in I}$ with respect to this restriction operation is formed as 
\[
\sum Y_i|a_i:=\left\{\sum y_i|a_i\colon y_i\in Y_i \text{ for each } i\in I\right\}
\]
for any partition $(a_i)_{i\in I}\in P(1)$.  

The nonempty set 
\[
\mathcal{P}(\mathbf{X}):=\{\mathbf{Y}|a\colon Y\in \mathcal{S}(\mathbf{X}), \; a\in \mathcal{A}\}. 
\]
is the carrier of a conditional set on $\mathcal{A}$, denoted by $\mathbf{P}(\mathbf{X})$, whose restriction operation is defined by 
\[
(\mathbf{Y}|a)|b:=\mathbf{Y}|a\wedge b, \quad b\in\mathcal{A}.  
\]
A concatenation inside $\mathbf{P}(\mathbf{X})$ is formed by 
\[
\sum_{i\in I} (\mathbf{Y}_i|b_i)|a_i:= \left(\mathbf{\sum_{i\in I} Y_i|a_i}\right)\Bigg|\vee_{i\in I} (a_i\wedge b_i) 
\]
where $\sum_{i\in I} Y_i|a_i$ is formed inside $\mathcal{S}(\mathbf{X})$; uniqueness is again guaranteed by \cite[Proposition 2.4]{drapeau2016algebra}.  
The conditional set $\mathbf{P(X)}$ is called the \emph{conditional power set} of $\mathbf{X}$.  
\end{definition}
\begin{remark}\label{superfluous}
The definition of a conditional set used in this note is different from the original definition in \cite{drapeau2016algebra}.  The only difference lies in dropping the property $x|a=y|b$ implies $a=b$, which is not needed. 

It is important to mention that the conditional power set constructed above is different from \cite[Definition 2.8]{drapeau2016algebra}.  
In \cite{drapeau2016algebra}, the conditional power set corresponds to the conditional set $\mathcal{S}(\mathbf{X})$ constructed above, and it was enough for the purposes of a conditional topology and conditional functional analysis developed in \cite{drapeau2016algebra}. However, for the development of a conditional measure and integration theory, it was realized that one needs to consider an extended conditional power set as is defined above, cf. \cite{jamneshan2018measures}.  It is also this extended conditional power set that coincides with the power set object in $\mathbf{Sh}(\mathcal{A},P)$.   
\end{remark}
It follows immediately from \cite[Theorem 2.9]{drapeau2016algebra} that the conditional power set has the structure of a complete Boolean algebra: 
\begin{theorem} \label{thm:powerset}
Let $\mathbf{X}$ be a conditional set. 
The classical inclusion relation defines a complete complemented distributive lattice structure on the carrier $\mathcal{P}(\mathbf{X})$ of the conditional power set of $\mathbf{X}$.  
\end{theorem}
\begin{definition}
Let $I$ be a nonempty index set, and let $\mathbf{X}_i$ be a conditional set on $\mathcal{A}_{a_i}$, $a_i\in\mathcal{A}$, for each $i\in I$.  
The conditional set consisting of objects $(x_i)_{i\in I}|a:=(x_i|a)_{i\in I}$ for $(x_i)_{i\in I}\in \prod_{i\in I} X_i$ and $a\in\mathcal{A}_{\wedge a_i}$ is called the \emph{conditional Cartesian product} of the family $(\mathbf{X}_i)$. 

Let $\mathbf{X}$ and $\mathbf{Y}$ be conditional sets on $\mathcal{A}_b$ for some $b\in \mathcal{A}$. A function $f$ from the carrier $X$ to the carrier $Y$ is said to be \emph{stable} if 
\[
     f\Big(\sum x_i|a_i\Big)=\sum f(x_i)|a_i, \; \forall (a_i)\in p(b), \; \forall (x_i) \subseteq X. 
\] 
Let $\mathbf{X}$ be a conditional set on $\mathcal{A}_a$ and $\mathbf{Y}$ be a conditional set on $\mathcal{A}_b$.  
A conditional subset $\mathbf{G_f}$ of  the conditional Cartesian product of $\mathbf{X}$ and $\mathbf{Y}$ is said to be the \emph{graph of a conditional function} $\mathbf{f}$ from $\mathbf{X}$ to $\mathbf{Y}$ if its carrier $G_f$ is the graph of a stable function of the underlying carriers (with respect to the relative algebra on which this product lives). 
\end{definition}
\begin{remark}
It can be verified that conditional sets form a category where the objects are conditional sets on relative algebras of a fixed complete Boolean algebra and the arrows are conditional functions.  
We call this category the \emph{category of conditional sets}, and denote it by $\textbf{cSet}(\mathcal{A})$. 
\end{remark}

\subsection{The equivalence of the categories of sheaves and conditional sets} \label{s13}

The proof of the following equivalence follows from the definitions. 
\begin{theorem}\label{equivcat}
\text{}
\begin{itemize}
\item[1.)]\label{shtocset} Let $\mathcal{X}$ be a sheaf over the site $(\mathcal{A},P)$. 
Recall that $a_\mathcal{X}$ is the largest element of $\mathcal{A}$ such that $X_{a_\mathcal{X}}$ is nonempty.  
For $x\in X_{a_\mathcal{X}}$ and $a\leq a_\mathcal{X}$, put
\[
x|a:= f^{a_\mathcal{X}}_a(x).
\]  
Then the collection 
\[
\mathbf{X}=\left\{x|a: x\in X_{a_\mathcal{X}}, \; a\in\mathcal{A}_{a_\mathcal{X}}\right\}
\]
is a conditional set on $\mathcal{A}_{a_\mathcal{X}}$. 
\item[2.)]\label{csettosh} Let $\mathbf{X}=\{x|b: x\in X, \; b\in \mathcal{A}_a\}$ be a conditional set on $\mathcal{A}_a$. Set 
\[
X_b=
\begin{cases}
\{x|b\colon x\in X\}, \quad &b\in\mathcal{A}_a, \\
\emptyset, &b\in \mathcal{A}\setminus \mathcal{A}_a. 
\end{cases}
\]
Let further $f^c_b \colon X_c \to X_b$ be defined by  
\[
\begin{cases}
x|c \to x|b, \quad & (b,c)\in\Delta_a, \\
\emptyset, \quad &(b,c)\in \Delta\setminus \Delta_a. 
\end{cases}
\]
Then $\mathcal{X}=((X_b)_{b\in\mathcal{A}}, (f^c_b)_{(b,c)\in\Delta})$ is a sheaf over the site $(\mathcal{A},P)$. 
\item[3.)] If $\mathcal{X}$ is a sheaf, $\mathbf{X}$ is the conditional set constructed from $\mathcal{X}$, and $\tilde{\mathcal{X}}$ is the sheaf constructed from $\mathbf{X}$, then $\mathcal{X}$ is isomorphic to $\tilde{\mathcal{X}}$. Conversely, if $\mathbf{X}$ is a conditional set, $\mathcal{X}$ is the sheaf constructed from $\mathbf{X}$, and $\tilde{\mathbf{X}}$ is the conditional set constructed from $\mathcal{X}$, then $\mathbf{X}$ is isomorphic to $\tilde{\mathbf{X}}$. 
\end{itemize}
Moreover, the power set objects and the Cartesian products in $\textbf{Sh}(\mathcal{A},P)$ correspond uniquely to the conditional power sets and conditional Cartesian products in $\textbf{cSet}$ respectively, and therefore $\textbf{Sh}(\mathcal{A},P)$ and $\textbf{cSet}$ are equivalent topoi.  
\end{theorem}

\subsection{Measure-theoretic sheaves}

We construct sheaves on the measure algebra associated to a $\sigma$-finite measure space.  For background on measure algebras, the interested reader is referred to \cite[Chapter 22]{monk1989handbook}.  
Let $(\Omega,\mathcal{F},\mu)$ be an arbitrary $\sigma$-finite measure space.  
Let $\mathcal{N}_\mu$ be the $\sigma$-ideal of $\mu$-negligible sets.  
The quotient Boolean algebra $\mathcal{A}=\mathcal{F}/\mathcal{N}_\mu$ is $\sigma$-complete and satisfies the countable chain condition and thus is complete. 
Partitions in $P(a)$ are at most countable for all $a\in\mathcal{A}$.  
We form the topos $\mathbf{Sh}(\mathcal{A},P)$ and name its objects \emph{measure-theoretic sheaves}.  
The naming is motivated by the instance that certain function spaces of measurable functions can be turned into measure-theoretic sheaves (similarly to sheaves over topological spaces).  
Notice that to any $\sigma$-subalgebra $\mathcal{G}\subseteq \mathcal{F}$, one could associate a topos $\mathbf{Sh}(\mathcal{B},P)$ where $\mathcal{B}=\mathcal{G}/(\mathcal{G}\cap \mathcal{N}_\mu)$.  
We also notice that any measure equivalent with respect to $\mu$ gives rise to the same topos.  
We provide an important class of examples of measure-theoretic sheaves. 
\begin{example}
Let $(M,\mathcal{M})$ be an arbitrary measurable space, and denote by $\mathscr{L}^0(M)$ the space of all measurable functions $x\colon \Omega\to M$. 
Let $L^0(M)$ denote the space of equivalence classes of $\mathscr{L}^0(M)$ modulo almost sure equality.  
Let $x|A$ denote the restriction of a function $x\in\mathscr{L}^0(M)$ to a measurable set $A\in\mathcal{F}$.  
The equivalence relation 
\[
x|A\sim y|B \quad \text{iff} \quad \mu(A\Delta B)=\mu(\{x\neq y\}\cap A)=0,
\]
where $\Delta$ denotes symmetric difference, is independent of representatives in $L^0(M)$. 
Therefore, one may  speak unambiguously of $x|a$ for $x\in L^0(M)$ and $a\in\mathcal{A}$.  
Let 
\[
\mathbf{M}:=\{x|a\colon x\in L^0(M),\; a\in\mathcal{A} \}. 
\]
For a countable partition $(A_n)$ of $\Omega$ consisting of measurable sets and a family $(x_n)$ in $\mathscr{L}^0(M)$ there exists a unique function $x\in  \mathscr{L}^0(M)$ such that $x|A_n=x_n|A_n$ for all $n$.  
This  shows that $\mathbf{M}$ is a conditional set on $\mathcal{A}$.  
With the help of  Theorem \ref{equivcat}, one can associate a measure-theoretic sheaf in $\mathbf{Sh}(\mathcal{A},P)$ to $\mathbf{M}$.  
\end{example}
In conditional set theory \cite{drapeau2016algebra}, a conditional version of basic results in analysis are proved by adapting classical definitions and proofs starting from classical elementary set theory. 
For example, to prove a conditional version of Tychonoff's theorem, one starts by defining a conditional topological space as a conditional set of conditional subsets closed under arbitrary conditional unions and conditionally finite conditional intersections, proceeds by defining conditional compactness of a conditional topological space using the conditional set operations, and finally verifies a conditional version of Tychonoff's theorem for these definitions, see \cite[Chapter 3]{drapeau2016algebra}, while adopting the classical definitions and arguments in a meaningful way.  
A conditional version of basic results in topology and functional analysis can be found in \cite{drapeau2016algebra} and a conditional version of basic results in measure theory can be found in \cite{jamneshan2018measures}.  
Next we want to apply the equivalence of topoi  in Theorem \ref{equivcat} to illustrate interpretations of some basic concepts of classical mathematics in the internal logic of $\mathbf{Sh}(\mathcal{A},P)$ by taking advantage of the worked out concepts in conditional set theory\footnote{In order to distinguish between classical concepts and classical theorems and their counterparts in the internal logic of $\mathbf{Sh}(\mathcal{A},P)$, we add the prefix $\mathbf{Sh}(\mathcal{A},P)$- to their names, for example $\mathbf{Sh}(\mathcal{A},P)$-compactness.}. 
\begin{enumerate}[i.]
\item The measure-theoretic sheaves associated to the function spaces\footnote{$\mathbb{N}$, $\mathbb{Z}$, $\mathbb{Q}$, $\mathbb{R}$ and $\mathbb{C}$ are endowed with their natural topologies.} $L^0(\mathbb{N})$, $L^0(\mathbb{Z})$, $L^0(\mathbb{Q})$, $L^0(\mathbb{R})$ and $L^0(\mathbb{C})$ are the $\mathbf{Sh}(\mathcal{A},P)$-natural numbers, -integers, -rational numbers, -real numbers and -complex numbers respectively.  
\item Any classical group that admits a conditional structure such that the group operation is a stable function induces a $\mathbf{Sh}(\mathcal{A},P)$-group. 
Conversely, any $\mathbf{Sh}(\mathcal{A},P)$-group can be interpreted as a classical group with the aforementioned property, cf. \cite[Definition 4.1]{drapeau2016algebra}.   
\item $\mathbf{Sh}(\mathcal{A},P)$-vector spaces over the $\mathbf{Sh}(\mathcal{A},P)$-real or complex numbers can be interpreted as classical $L^0(\mathbb{R})$- or  $L^0(\mathbb{C})$-modules satisfying a stability property, cf. \cite[Section 5.1]{drapeau2016algebra}.  
\item Let $\mathbf{X}$ be a conditional set.  A stable collection of stable sets in $\mathcal{S}(\mathbf{X})$ which is a classical topological base induces a $\mathbf{Sh}(\mathcal{A},P)$-topological base on the related sheaf, cf. \cite[Proposition 3.5]{drapeau2016algebra}.  
\item Let $\mathbf{X}$ be a conditional set.  A stable collection of stable sets in $\mathcal{S}(\mathbf{X})$ which is a classical filter base induces a $\mathbf{Sh}(\mathcal{A},P)$-filter base on the related sheaf, cf. \cite[Proposition 3.15]{drapeau2016algebra}. 
\item By the previous two items and by \cite[Propositions 3.11 and 3.22]{drapeau2016algebra}, we can establish equivalences between classical continuity and convergence and $\mathbf{Sh}(\mathcal{A},P)$-continuity and convergence.  This also leads to a classical interpretations of $\mathbf{Sh}(\mathcal{A},P)$-compactness (cf. \cite[Section 3.4]{drapeau2016algebra}).  
\item Let $E$ be an $L^0(\mathbb{R})$-module that admits a conditional structure\footnote{The measure algebra acts on $E$ by $L^0(\mathbb{R})$-scalar multiplication $1_A\cdot x$ where $A$ is a representative of $a\in\mathcal{A}$ and $x\in E$. We say that $E$ admits the structure of a conditional set if for every partition $(a_n)\in P(1)$ and $(x_n)\subseteq E$ there exists a unique $x\in E$ such that $1_{A_n} x=1_{A_n}x_n$ for all $n$. In this case, define on $E\times \mathcal{A}$ the equivalence relation $(x,a)\sim (y,b)$ if $a=b$ and $1_Ax=1_By$, and denote the resulting equivalence classes by $x|a$. Then $\mathbf{E}:=\{x|a\colon x\in E, a\in\mathcal{A}\}$ has the structure of a conditional vector space.}  and let $n\colon E\to L^0(\mathbb{R}_+)$ be a vector norm.  It can be shown that the function $n$ is stable, and thus induces a conditional norm $\mathbf{n}\colon \mathbf{E}\to \mathbf{R_+}$.  
The collection 
\[
B_r(x):=\{y\in E\colon n(x-y)<r\}, \quad r\in L^0(\mathbb{R}_{++}),\; x\in E,
\]
is in $\mathcal{S}(\mathbf{E})$ and defines a classical topological base and corresponds to the conditional topology on $\mathbf{E}$ induced by the conditional norm.  The  associated sheaf is a $\mathbf{Sh}(\mathcal{A},P)$-normed space. 

For example, let $M$ be a Banach space and $L^0(M)$ be the space of strongly measurable functions $x\colon \Omega\to M$. Then the norm of $M$ extends to a vector norm on $L^0(M)$, and the measure-theoretic sheaf associated to $L^0(M)$ is a $\mathbf{Sh}(\mathcal{A},P)$-normed space. Actually, it can be shown that $L^0(M)$ is the classical interpretation of a $\mathbf{Sh}(\mathcal{A},P)$-Banach space.  
\item  A $\mathbf{Sh}(\mathcal{A},P)$-complex Hilbert space can be interpreted as a classical complete topological $L^0(\mathbb{C})$-module $E$ satisfying the countable gluing property endowed with  an $L^0(\mathbb{C})$-valued inner product that induces the topology. 
 Conversely, the sheaf associated to such topological module is a $\mathbf{Sh}(\mathcal{A},P)$-complex Hilbert space.  
 For any classical complex Hilbert space $H$, the sheaf associated to the function space $L^0(H)$ of equivalence classes of strongly measurable functions with values in $H$ is a $\mathbf{Sh}(\mathcal{A},P)$-complex Hilbert space. 

 Let $\mathcal{G}$ be a sub-$\sigma$-algebra of $\mathcal{F}$.  
For a nonnegative random variable $x\colon \Omega\to \mathbb{R}$, the conditional expectation operator $\mathbf{E}(x\,|\,\mathcal{G})$ is well defined with values in $[0,\infty]$.  
We denote by $L^0(\mathcal{G},\mathbb{C})$ the space of equivalence classes of $\mathcal{G}$-measurable random variables $x\colon \Omega\to \mathbb{C}$.  
Define 
\[
L^2_{\mathcal{F}|\mathcal{G}}:=\{x\in L^0(\mathcal{F},\mathbb{C})\colon \mathbf{E}(|x|^2 \,|\, \mathcal{G})\in L^0(\mathcal{G},\mathbb{R})\}
\]
Then $L^2_{\mathcal{F}|\mathcal{G}}$ is an $L^0(\mathbb{C})$-module closed under countable gluings endowed with an $L^0(\mathcal{G},\mathbb{C})$-valued inner product $\mathbf{E}\left(x\bar{y}\,|\,\mathcal{G}\right)$. 
The associated sheaf is a Hilbert space in $\mathbf{Sh}(\mathcal{B},P)$ where $\mathcal{B}=\mathcal{G}/(\mathcal{G}\cap \mathcal{N}_\mu)$.  
 \end{enumerate}

\end{document}